\documentclass[a4paper,11pt]{amsart}
\usepackage{amssymb}
\usepackage{mathrsfs}
\usepackage{amsmath,amssymb,amsthm,latexsym,amscd,mathrsfs}
\usepackage{indentfirst}
\usepackage{stmaryrd}
 \newcommand{\ROM}[1]{\mathrm{\uppercase\expandafter{\romannumeral#1}}}
  \theoremstyle{definition}

 \newtheorem{thm}{Theorem}[section]

 \newtheorem{rem}{Remark}[section]

\newtheorem{ack}{Acknowledgements}   

\title[Isoparametric foliation and Yau conjecture on the first eigenvalue, II]{\textbf{Isoparametric foliation and Yau conjecture on the first eigenvalue, II}}
\author[Z.Z.Tang]{Zizhou Tang}\address{School of Mathematical Sciences, Laboratory of Mathematics and Complex Systems, Beijing Normal
University, Beijing 100875, China}\email{zztang@bnu.edu.cn}
\thanks {The project is partially supported by the NSFC (No.11071018,
No.11331002 and No.11301027), the SRFDP (No.20130003120008).}
\author[Y.Q.Xie]{Yuquan Xie}
\address{Department of Mathematics, Hangzhou Normal University, Zhejiang 310036,
China} \email{yuqxie@pku.edu.cn}
\author[W. J. Yan]{Wenjiao Yan*}
\thanks {* The third author is the corresponding author.}
\address{School of Mathematical Sciences, Laboratory of Mathematics and Complex Systems, Beijing Normal
University, Beijing 100875, China} \email{wjyan@bnu.edu.cn}

\date{}
\begin{document}

\maketitle
\begin{abstract}
This is a continuation of \cite{TY}, which investigated the
first eigenvalues of minimal isoparametric hypersurfaces with $g=4$
distinct principal curvatures and focal submanifolds in unit
spheres. For the focal submanifolds with $g=6$, the present paper
obtains estimates on all the eigenvalues, among others, giving an
affirmative answer in one case to the problem posed in \cite{TY},
which may be regarded as a generalization of Yau's conjecture. In
two of the four unsettled cases in \cite{TY} for focal submanifolds
$M_1$ of OT-FKM-type, we prove the first eigenvalues to be their
dimensions, respectively.
\end{abstract}

%%%%%%%%%%%%%%%%%%%%%%%%%%%%%%%%%%%%%%%%%%%%%%%%%%%%%%%%%%%%%%%%%%%%%%%%%%%%%%%%%%%%%%%%%%%%%%%%%%%%%%%%%%%%%%%%%%%%%%%%%%%%%%%%%%%%%%%%%%%%%
\section{\textbf{Introduction}}
Let $M^n$ be an $n$-dimensional closed connected Riemannian manifold
and $\Delta $ be the Laplace-Beltrami operator acting on a
$C^{\infty}$ function $f$ on $M$ by $\Delta f$ $=-$ div$(\nabla f)$,
the negative of divergence of the gradient $\nabla f$. It is well
known that $\Delta$ is an elliptic operator and has a discrete
spectrum
$$\{0=\lambda_0(M)<\lambda_1(M)\leqslant \lambda_2(M)\leqslant\cdots \leqslant\lambda_k(M),\cdots, \uparrow \infty\}$$
with each eigenvalue occurs as many times as its multiplicity.
As usual, we call $\lambda_1(M)$ the first eigenvalue of $M$. A well known conjecture of
S.T.Yau states that
\vspace{2mm}

\noindent \textbf{Yau conjecture (\cite{Yau}):}\,\, {\itshape The
first eigenvalue of every (embedded) closed minimal hypersurface
$M^n$ in the unit sphere $S^{n+1}(1)$ is just $n$. } \vspace{2mm}

Up to now, Yau's conjecture is still far from being solved. The most
recent contribution to this problem is given by \cite{TY}. They give
an affirmative answer to Yau's conjecture for closed minimal
isoparametric hypersurfaces $M^n$ in $S^{n+1}(1)$.

By definition, an isoparametric hypersurface $M^n$ in the unit
sphere $S^{n+1}(1)$ is a hypersurface with constant principal
curvatures. We denote by $g$ the number of distinct principal
curvatures, and $m_1, m_2$ their multiplicities (details will be discussed in the next
section).

In fact, for the minimal isoparametric hypersurfaces with $g=6$, the
proof of Yau's conjecture is just a simple combination of the
results of \cite{MOU}, \cite{Kot} with the classification theorems
of \cite{DN} and \cite{Miy2}, \cite{Miy}. An interesting problem
naturally arises as to whether it is possible to give a direct proof
without using the classification theorems of Dorfmeister-Neher and
Miyaoka, which states that all the isoparametric hypersurfaces with
$g=6$ in unit spheres are homogeneous. As the first result of this
paper, we provide a direct proof. Moreover, we obtain more
information than that in \cite{MOU}, which only focused on the first
eigenvalue of the minimal homogeneous hypersurfaces.

 \vspace{2mm}

\begin{thm}\label{prop 1 of TXY}
Let $M^{12}$ be a closed minimal isoparametric hypersurface in
$S^{13}(1)$ with $g=6$ and $(m_1, m_2)=(2, 2)$. Then
$$\lambda_1(M^{12})=12$$
with multiplicity $14$. Furthermore, we have the inequality
$$\lambda_k(M^{12})>\frac{3}{7}~\lambda_k(S^{13}(1)), \qquad k=1, 2, \cdots.$$
\end{thm}
\vspace{2mm}

Other than the minimal isoparametric hypersurfaces, \cite{TY}
originally studied the first eigenvalues of the focal submanifolds
of the isoparametric foliation in $S^{n+1}(1)$, which are in fact
the minimal submanifolds in $S^{n+1}(1)$. \vspace{2mm}

\noindent \textbf{\emph{Theorem 1.3 in \cite{TY}}.}\,\, \emph{Let $M_1$ be
the focal submanifold of an isoparametric hypersurface with $g=4$ in
$S^{n+1}(1)$. If $\dim M_1\geqslant \frac{2}{3}n+1$, then
$$\lambda_1(M_1)=\dim M_1$$ with multiplicity $n+2$. A similar conclusion holds for the other focal submanifold $M_2$.}
\vspace{2mm}

As asserted in \cite{TY}, there are only four unsettled cases for
the first eigenvalues of the focal submanifolds $M_1$ (\emph{i.e.},
$f^{-1}(1)$, $f$ is the restriction of the OT-FKM polynomial on the
unit sphere) in the isoparametric foliation of OT-FKM type ($g=4$).
Namely, $(m_1, m_2)=(1, 1), (4, 3)$ associated with one homogeneous
and one inhomogeneous examples, and $(5, 2)$. Unfortunately, their
method is invalid for these cases. As the next aim of this paper, we
consider $M_1$ with multiplicity pairs $(m_1, m_2)=(1, 1)$, or $(4,
3)$ associated with the homogeneous example, to obtain one of our
main results as follows. \vspace{2mm}

\begin{thm}\label{thm2 TXY focal g=4}
For the focal submanifold $M_1$ of OT-FKM type in $S^5(1)$ with
$(m_1, m_2)=(1, 1)$,
$$\lambda_1(M_1)=\dim M_1=3$$
with multiplicity $6$; for the focal submanifold $M_1$ of
homogeneous OT-FKM type in $S^{15}(1)$ with $(m_1, m_2)=(4, 3)$,
$$\lambda_1(M_1)=\dim M_1=10$$
with multiplicity $16$.
\end{thm}
\vspace{2mm}

\begin{rem}
As asserted in \cite{TY}, the first eigenvalue of the focal
submanifold $M_2$ of OT-FKM type in $S^5(1)$ with $(m_1, m_2)=(1,
1)$ is equal to its dimension. As for the focal submanifold $M_2$ of
homogeneous OT-FKM type in $S^{15}(1)$ with $(m_1, m_2)=(4, 3)$, its
dimension satisfies the assumption of Theorem 1.3 in \cite{TY}; thus
the first eigenvalue is equal to its dimension. By virtue of
eigenfunctions constructed by Solomon on $M_1$ of OT-FKM type with
$(m_1, m_2)=(5, 2)$, we see that the first eigenvalue is less than
its dimension (\emph{cf}. \cite{Sol}).
\end{rem}

Notice that in their method calculating the first eigenvalues of the
focal submanifolds, \cite{TY} took average value of the gradient of
the test functions at each pair of antipodal points. However, in the
case $g=6$, the average value is not accurate enough to meet our
requirement. In this paper, by investigating the shape operators of
the focal submanifolds, we obtain estimates on the first
eigenvalues.
 \vspace{2mm}

\begin{thm}\label{thm1 TXY focal g=6}
For the focal submanifolds of an isoparametric foliation with $g=6$,
we have
\begin{itemize}

\item[(i)] when $(m_1, m_2)=(1, 1)$, the first eigenvalues of the focal submanifolds
$M_1$ and $M_2$ in $S^{7}(1)$ satisfy
$$3\leqslant\lambda_1(M_1), \,\,\, \lambda_1(M_2)\leqslant \dim M_1=\dim M_2=5.$$

\item[(ii)] when $(m_1, m_2)=(2, 2)$, the $k$-th eigenvalues of the focal submanifolds
$M_1$ and $M_2$ in $S^{13}(1)$ satisfy
$$\lambda_k(S^{13}(1))\leqslant(3+\frac{99\sqrt{3}}{40\pi})\cdot\lambda_k(M_1),\quad \lambda_k(S^{13}(1))\leqslant (6-\frac{117\sqrt{3}}{20\pi})\cdot\lambda_k(M_2),$$
for $k=1, 2, \cdots.$ In particular,
\begin{equation}\label{g=6 m=2 focal}
\lambda_1(M_2)=\dim M_2=10
\end{equation} with multiplicity $14$.
\end{itemize}

\end{thm}

\begin{rem}
In the case $g=6$ and $(m_1, m_2)=(2, 2)$, we will distinguish $M_1$
from $M_2$ in Section 4, following the notations in \cite{Miy}. The
equality (\ref{g=6 m=2 focal}) in Theorem \ref{thm1 TXY focal g=6}
gives in this case an affirmative answer to the problem in
\cite{TY}, which may be regarded as a generalization of Yau's
conjecture. Unfortunately, we have not obtained the accurate value
of the first eigenvalue of $M_1$. Notice that $M_1$ and $M_2$ are
not congruent in $S^{13}(1)$ (\emph{cf}. \cite{Miy}). In fact,
comparing the Ricci tensors by Gauss equation, one finds that $M_1$
and $M_2$ are not isometric. The problem of determining
$\lambda_1(M_1)$ is still open!
\end{rem}

\section{\textbf{Preliminary}}

An oriented hypersurface $M^n$ in the unit sphere $S^{n+1}(1)$ with
constant principal curvatures is called an isoparametric
hypersurface (\emph{cf.} \cite{Car1}, \cite{Car2}, \cite{CR}). It is
well known that a closed isoparametric hypersurface is an oriented,
embedded hypersurface. Denote by $\xi$ a unit normal vector field
along $M^n$ in $S^{n+1}(1)$, $g$ the number of distinct principal
curvatures of $M$, $\cot \theta_{\alpha}~ (\alpha=1,...,g;~
0<\theta_1<\cdots<\theta_{g} <\pi)$ the principal curvatures with
respect to $\xi$ and $m_{\alpha}$ the multiplicity of $\cot
\theta_{\alpha}$. According to M\"{u}nzner (\cite{Mun}), the number
$g$ must be $1, 2, 3, 4$ or $6$; $m_{\alpha}=m_{\alpha+2}$ (indices
mod $g$) and $\theta_{\alpha}=\theta_1+\frac{\alpha-1}{g}\pi$
$(\alpha = 1,..., g)$.

For isoparametric hypersurfaces in unit spheres with $g = 1, 2, 3$,
Cartan classified them to be homogeneous (\emph{cf. }\cite{Car1},
\cite{Car2}); when $g = 6$, Abresch (\cite{Abr}) showed that the
multiplicity of each principal curvatures only takes values
$m_1=m_2=1$ or $2$. Dorfmeister-Neher (\cite{DN}) and Miyaoka
(\cite{Miy}) proved the homogeneity of such hypersurfaces,
respectively; for the most complicated case $g=4$, Cecil-Chi-Jensen
(\cite{CCJ}), Immervoll (\cite{Imm}) and Chi (\cite{Chi}) proved a
far reaching result that they are either homogeneous or of
OT-FKM-type except possibly for the case $(m_1, m_2)=(7,8)$.

A well known result of Cartan states that isoparametric hypersurfaces come as a family of
parallel hypersurfaces. To be more specific, given an isoparametric hypersurface $M^n$ in $S^{n+1}(1)$ and a smooth field $\xi$ of unit
normals to $M$, for each $x\in M$ and $\theta\in \mathbb{R}$, we can define
$\phi_{\theta}: M^n\rightarrow S^{n+1}(1)$ by
$$\phi_{\theta}(x)=\cos \theta~ x +\sin \theta~ \xi(x).$$
Clearly, $\phi_{\theta}(x)$ is the point at an oriented distance
$\theta$ to $M$ along the normal geodesic through $x$. If
$\theta\neq \theta_{\alpha}$ for any $\alpha=1,...,g$,
$\phi_{\theta}$ is a parallel hypersurface to $M$ at an oriented
distance $\theta$, which we will denote by $M_{\theta}$
henceforward. If $\theta= \theta_{\alpha}$ for some
$\alpha=1,...,g$, it is easy to find that for any vector $X$ in the
principal distributions $E_{\alpha}(x)=\{X\in T_xM ~|~A_{\xi}X=\cot
\theta_{\alpha}X\}$, where $A_{\xi}$ is the shape operator with
respect to $\xi$, $(\phi_{\theta})_{\ast}X=(\cos \theta-\sin \theta
\cot
\theta_{\alpha})X=\frac{sin(\theta_{\alpha}-\theta)}{sin\theta_{\alpha}}
X=0$. In other words, in case that $\cot \theta=\cot
\theta_{\alpha}$ is a principal curvature of $M$, $\phi_{\theta}$ is
not an immersion, whose image is actually a \emph{focal submanifold}
of codimension $m_{\alpha}+1$ in $S^{n+1}(1)$.

As asserted by M\"{u}nzner, regardless of the number of distinct
principal curvatures of $M$, there are only two distinct focal
submanifolds in a parallel family of isoparametric hypersurfaces,
and every isoparametric hypersurface is a tube of constant radius
over each focal submanifold. Denote by $M_1$ the focal submanifold
in $S^{n+1}(1)$ at an oriented distance $\theta_1$ along $\xi$ from
$M$ with codimension $m_1+1$, $M_2$ the focal submanifold in
$S^{n+1}(1)$ at an oriented distance $\frac{\pi}{g}-\theta_1$ along
$-\xi$ from $M$ with codimension $m_2+1$. Another choice of the
normal direction will lead to the exchange between the focal
submanifolds $M_1$ and $M_2$. In virtue of Cartan's identity, one
sees that both the focal submanifolds $M_1$ and $M_2$ are minimal in
$S^{n+1}(1)$ (\emph{cf.} \cite{CR}).

\vspace{2mm}

\section{\textbf{Isoparametric hypersurfaces with $(g, m_1, m_2)=(6, 2, 2)$.}}

Let $\phi: M^n \rightarrow S^{n+1}(1)(\subset \mathbb{R}^{n+2})$ be
a closed isoparametric hypersurface and again $M_{\theta}$ be the
parallel hypersurface defined by $\phi_{\theta}: M^n \rightarrow
S^{n+1}(1)~(-\pi<\theta<\pi, \cot\theta\neq \cot \theta_{\alpha})$,
$$\phi_{\theta}(x) = \cos \theta~ x + \sin \theta~ \xi(x).$$
It is clear that for $X\in E_{\alpha}$,
\begin{equation}\label{phi theta star}
(\phi_{\theta})_{\ast} X=\frac{\sin(\theta_{\alpha}-\theta)}{\sin
\theta_{\alpha}} \widetilde{X},
\end{equation}
where $\widetilde{X} \sslash X$ as vectors in $\mathbb{R}^{n+2}$.

Following \cite{TY}, we will apply the theorem below to the case
$V=S^{n+1}(1)$ and $W=M_1\cup M_2$ and prove Theorem \ref{prop 1 of
TXY} by estimating the eigenvalue $\lambda_{k}(M^n)$ from below.
\vspace{1mm}

\noindent \textbf{Theorem (Chavel and Feldman \cite{CF}, Ozawa
\cite{Oza})}\,\, {\itshape Let $V$ be a closed, connected Riemannian
manifold and $W$ a closed submanifold. For any sufficiently small
$\varepsilon>0$, set $W(\varepsilon)=\{x\in V:~ dist(x,
W)<\varepsilon\}$. Let $\lambda^D_k(\varepsilon)$ $(k=1, 2, ...)$ be
the $k$-th eigenvalue on $V-W(\varepsilon)$ under the Dirichlet
boundary condition. If $\dim V\geqslant \dim W+ 2$, then
\begin{equation}\label{lambda k limit}
\lim_{\varepsilon\to 0}\lambda^D_k(\varepsilon)=\lambda_{k-1}(V).
\end{equation}}

\noindent \textbf{\emph{Proof of Theorem \ref{prop 1 of TXY}}}. In
our case with $(g, m_1, m_2)=(6, 2, 2)$, denote by $M^{12}$ the
minimal isoparametric hypersurface. Clearly,
$\theta_1=\frac{\pi}{12}$. For sufficiently small $\varepsilon> 0$,
set
$$M(\varepsilon)=\bigcup_{\theta\in[-\frac{\pi}{12}+\varepsilon, ~\frac{\pi}{12}-\varepsilon]}M_{\theta},$$
which is a tube around $M^{12}$. According to the previous theorem,
\begin{equation}\label{Chavel}
\lim_{\varepsilon\to
0}\lambda^D_{k+1}(M(\varepsilon))=\lambda_{k}(S^{13}(1)), \qquad
k=1, 2, \cdots.
\end{equation}

Let $\Big\{\widetilde{e}_{\alpha,i}~ \boldsymbol{|} ~ i=1, 2,~
\alpha=1,..,6,~ \widetilde{e}_{\alpha, i}\in E_{\alpha}\Big\}$ be a
local orthonormal frame field on $M$. Then
$$\Big\{\frac{\partial}{\partial\theta},~
e_{\alpha,i}~\boldsymbol{|}~
e_{\alpha,i}=\frac{\sin\theta_{\alpha}}{\sin(\theta_{\alpha}-\theta)}\widetilde{e}_{\alpha,i},~
i=1, 2,~ \alpha=1,..,6,~ \theta\in[-\frac{\pi}{12}+\varepsilon,
~\frac{\pi}{12}-\varepsilon]\Big\}$$ constitutes a local orthonormal
frame field on $M(\varepsilon)$. From the formula (\ref{phi theta
star}), we derive the following equality up to a sign:
\begin{equation}\label{volume}
dM(\varepsilon)=16\cos^22\theta\sin^2(\frac{\pi}{6}+2\theta)\sin^2(\frac{\pi}{6}-2\theta)d\theta dM
\end{equation}
where $dM(\varepsilon)$ and $dM$ are the volume elements of $M(\varepsilon)$ and $M$, respectively.

Again following \cite{TY}, let $h$ be a nonnegative, increasing
smooth function on $[0, \infty)$ satisfying $h=1$ on $[2, \infty)$
and $h=0$ on $[0, 1]$. For sufficiently small $\eta>0$, let
$\psi_{\eta}$ be a nonnegative smooth function on $[\eta,
\frac{\pi}{2}-\eta]$ such that

$(i)$ $\psi_{\eta}(\eta)=\psi_{\eta}(\frac{\pi}{2}-\eta)=0$,

$(ii)$ $\psi_{\eta}$ is symmetric with respect to $x=\frac{\pi}{4}$,

$(iii)$ $\psi_{\eta}(x)=h(\frac{x}{\eta})$ on $[\eta, \frac{\pi}{4}]$.

\noindent
Let $f_k$ $(k=0, 1, ...)$ be the $k$-th eigenfunctions on $M$ which are orthogonal to
each other with respect to the square integral inner product on $M$ and $L_{k+1} =Span\{f_0, f_1,..., f_k\}$.

For each fixed $\theta\in[-\frac{\pi}{12}+\varepsilon,
~\frac{\pi}{12}-\varepsilon]$, denote
$\pi=\pi_{\theta}=\phi^{-1}_{\theta}: M_{\theta}\rightarrow M$. Then
any $\varphi\in L_{k+1}$ on $M$ can give rise to a function
$\Phi_{\varepsilon} : M(\varepsilon)\rightarrow \mathbb{R}$ by
$$\Phi_{\varepsilon}(x) = \psi_{3\varepsilon}(3(\frac{\pi}{12}-\theta))(\varphi\circ\pi)(x),$$
where $\theta$ is characterized by $x\in M_{\theta}$. It is easily
seen that $\Phi_{\varepsilon}$ is a smooth function on
$M(\varepsilon)$ satisfying the Dirichlet boundary condition and
square integrable.

By the mini-max principle, we obtain:
\begin{equation}\label{min max}
\lambda^D_{k+1}(M(\varepsilon))\leqslant \sup_{\varphi\in L_{k+1}}\frac{\|\nabla \Phi_{\varepsilon}\|_2^2}{\|\Phi_{\varepsilon}\|_2^2}.
\end{equation}
In the following, we will concentrate on the calculation of
$\frac{\|\nabla
\Phi_{\varepsilon}\|_2^2}{\|\Phi_{\varepsilon}\|_2^2}$. Observing
that the normal geodesic starting from $M$ is perpendicular to each
$M_{\theta}$, we obtain
$$\|\nabla \Phi_{\varepsilon}\|_2^2=\int_{M(\varepsilon)}9(\psi^{\prime}_{3\varepsilon})^2\varphi(\pi)^2dM(\varepsilon)
+\int_{M(\varepsilon)}\psi^{2}_{3\varepsilon}|\nabla \varphi(\pi)|^2dM(\varepsilon).$$
On the other hand, a simple calculation leads to
\begin{eqnarray}\label{Phi length}
\|\Phi_{\varepsilon}\|_2^2&=&\int_{M(\varepsilon)}\psi^2_{3\varepsilon}(3(\frac{\pi}{12}-\theta))\varphi(\pi(x))^2 dM(\varepsilon)\nonumber\\
&=&\int_M\int_{-\frac{\pi}{12}+\varepsilon}^{\frac{\pi}{12}-\varepsilon}16\cos^22\theta\sin^2(\frac{\pi}{6}+2\theta)
\sin^2(\frac{\pi}{6}-2\theta)\psi^2_{3\varepsilon}(3(\frac{\pi}{12}-\theta))\varphi(\pi(x))^2d\theta dM\nonumber\\
&=&\frac{16}{3}\|\varphi\|^2_2\Big(\int_{3\varepsilon}^{\frac{\pi}{2}-3\varepsilon}\psi^2_{3\varepsilon}(x)\sin^{2}(\frac{2}{3}x)
\cos^{2}(\frac{\pi}{6}-\frac{2}{3}x)\sin^2(\frac{\pi}{3}-\frac{2}{3}x)~dx\Big).\nonumber
\end{eqnarray}
\noindent
For the sake of convenience, let us decompose
\begin{equation}\label{Phi}
\frac{\|\nabla \Phi_{\varepsilon}\|_2^2}{\|\Phi_{\varepsilon}\|_2^2}=I(\varepsilon)+II(\varepsilon),
\end{equation}
with
\begin{eqnarray}\label{I}
I(\varepsilon)
&=& \frac{\int_{M(\varepsilon)}9(\psi^{\prime}_{3\varepsilon})^2\varphi(\pi)^2dM(\varepsilon)}{\int_{M(\varepsilon)}(\psi_{3\varepsilon})^2\varphi(\pi)^2~dM(\varepsilon)}\\
&=&\frac{9\int_{3\varepsilon}^{\frac{\pi}{2}-3\varepsilon}(\psi^{\prime}_{3\varepsilon}(x))^2\sin^{2}(\frac{2}{3}x)
\cos^{2}(\frac{\pi}{6}-\frac{2}{3}x)\sin^2(\frac{\pi}{3}-\frac{2}{3}x)~dx}
{\int_{3\varepsilon}^{\frac{\pi}{2}-3\varepsilon}\psi^2_{3\varepsilon}(x)\sin^{2}(\frac{2}{3}x)
\cos^{2}(\frac{\pi}{6}-\frac{2}{3}x)\sin^2(\frac{\pi}{3}-\frac{2}{3}x)~ dx}\nonumber
\end{eqnarray}
and
\begin{equation}\label{II}
II(\varepsilon)
 =\frac{\int_{M(\varepsilon)}\psi^{2}_{3\varepsilon}|\nabla \varphi(\pi)|^2dM(\varepsilon)}
 {\int_{M(\varepsilon)}\psi^2_{3\varepsilon}\varphi(\pi)^2 dM(\varepsilon)}.
%&=&\frac{\|\nabla \varphi(\pi)\|^2_2}{\|\varphi(\pi)\|_2^2} \frac{4\int_{2\varepsilon}^{\frac{\pi}{2}-2\varepsilon}(\psi^{\prime}_{2\varepsilon}(x))^2\sin^{m_1}x\cos^{m_2}x dx}{\int_{2\varepsilon}^{\frac{\pi}{2}-2\varepsilon}\psi^2_{2\varepsilon}(x)\sin^{m_1}x\cos^{m_2}x dx}\nonumber\\
\end{equation}

Firstly, as in \cite{TY}, we deduce without difficulty that
\begin{equation}\label{limit I}
\lim_{\varepsilon\rightarrow 0}I(\varepsilon)=0.
\end{equation}

Next, we turn to the estimate on $II(\varepsilon)$. Decompose
$\nabla \varphi = Z_1 + \cdots + Z_6 \in E_1\oplus \cdots\oplus
E_6$, and set $k_{\alpha}=\frac{sin(\theta_{\alpha}-\theta)}{sin
\theta_{\alpha}}$ for $\alpha=1,...,6$. It follows obviously that
\begin{equation}\label{nabla varphi(pi)}
\left\{ \begin{aligned}
\quad|\nabla \varphi|^2~~~&=|Z_1|^2+\cdots+|Z_6|^2\\
|\nabla
\varphi(\pi)|^2&=\frac{1}{k_1^2}|Z_1|^2+\cdots+\frac{1}{k_6^2}|Z_6|^2.
\end{aligned}\right.
\end{equation}
Moreover, for $\alpha=1,...,6$, define
\begin{eqnarray}\label{K alpha}
K_{\alpha}&:=&16\int_{-\frac{\pi}{12}}^{\frac{\pi}{12}}\frac{\cos^{2}(2\theta)
\sin^{2}(\frac{\pi}{6}+2\theta)\sin^2(\frac{\pi}{6}-2\theta)}{k^2_{\alpha}}d\theta\\
%&=&\sin^2\theta_{\alpha}\int_0^{\frac{\pi}{4}}\frac{\sin^{m_1}2x\cos^{m_2}2x}{\sin^2(\frac{\alpha-1}{4}\pi+x)} dx,\nonumber\\
G&:=&32\int_{-\frac{\pi}{12}}^{\frac{\pi}{12}}\cos^{2}(2\theta)
\sin^{2}(\frac{\pi}{6}+2\theta)\sin^2(\frac{\pi}{6}-2\theta)~d\theta.\nonumber
%&=&2\int_0^{\frac{\pi}{4}}\sin^{m_1}2x\cos^{m_2}2x~dx.\nonumber
\end{eqnarray}
Let $\displaystyle K=\max_{\alpha}\{K_{\alpha}\}$. Then combining
(\ref{Phi}), (\ref{I}), (\ref{II}), (\ref{limit I}), (\ref{nabla
varphi(pi)}) with (\ref{K alpha}), we accomplish that
\begin{equation}\label{limit II}
\lim_{\varepsilon\rightarrow 0}\frac{\|\nabla
\Phi_{\varepsilon}\|_2^2}{\|\Phi_{\varepsilon}\|_2^2}= \frac{
\sum_{\alpha}K_{\alpha}\|Z_{\alpha}\|_2^2}{\|\varphi\|_2^2\cdot
\frac{1}{2}G} \leqslant \frac{2K}{G}\cdot \frac{\|\nabla
\varphi\|_2^2}{\|\varphi\|_2^2}.
\end{equation}

\noindent
Therefore, putting (\ref{Chavel}), (\ref{min max}) and (\ref{limit II}) together, we obtain
\begin{equation}\label{sphere M}
\lambda_k(S^{13}(1))
=\lim_{\varepsilon\rightarrow 0}\lambda_{k+1}^D(M(\varepsilon))
\leqslant\lim_{\varepsilon\rightarrow 0}\sup_{\varphi\in L_{k+1}}\frac{\|\nabla \Phi_{\varepsilon}\|_2^2}{\|\Phi_{\varepsilon}\|_2^2}
\leqslant\lambda_k(M^{12}) \frac{2K}{G}.
\end{equation}

\noindent Comparing the leftmost side with the rightmost side of
(\ref{sphere M}), we find a sufficient condition to complete the
proof of Theorem \ref{prop 1 of TXY}, namely,
\begin{equation}\label{KG}
K<\frac{7}{6}G,
\end{equation}
since then $\lambda_{15}(S^{13}(1))=28<
\lambda_{15}(M^{12})\cdot\frac{7}{3}$, which implies immediately
that $\lambda_{15}(M^{12})> 12$. On the other hand, recall that $12$
is an eigenvalue of $M^{12}$ with multiplicity at least $14$.
Therefore, the first eigenvalue of $M^{12}$ must be $12$ with
multiplicity $14$.

We are left to verify the inequality (\ref{KG}). Observing that
$K_1=K_6$, $K_2=K_5$ and $K_3=K_4$, we give the following
straightforward verification.

$(i)$ \begin{eqnarray}\label{K1}
K_1&=&16\int_{-\frac{\pi}{12}}^{\frac{\pi}{12}}\frac{\cos^{2}(2\theta)
\sin^{2}(\frac{\pi}{6}+2\theta)\sin^2(\frac{\pi}{6}-2\theta)\sin^2{\frac{\pi}{12}}}{\sin^2(\frac{\pi}{12}-\theta)}~d\theta\nonumber\\
%&=& 16(2-\sqrt{3})\int_{-\frac{\pi}{12}}^{\frac{\pi}{12}}\cos^{2}(2\theta)
%\sin^{2}(\frac{\pi}{6}+2\theta)\cos^2(\frac{\pi}{12}-\theta)~d\theta\nonumber\\
&=& 16(2-\sqrt{3})(\frac{\pi}{64}+\frac{63\sqrt{3}}{1280}),\nonumber
\end{eqnarray} while
\begin{equation*}\label{G}
G~~=32\int_{-\frac{\pi}{12}}^{\frac{\pi}{12}}\cos^{2}(2\theta)
\sin^{2}(\frac{\pi}{6}+2\theta)\sin^2(\frac{\pi}{6}-2\theta)~d\theta=\frac{\pi}{6}.\qquad\quad
\end{equation*}
Therefore, $$K_1<\frac{7}{6}G.$$

$(ii)$ \begin{eqnarray*}
  K_2 &=& 16\int_{-\frac{\pi}{12}}^{\frac{\pi}{12}}\frac{\cos^{2}(2\theta)
\sin^{2}(\frac{\pi}{6}+2\theta)\sin^2(\frac{\pi}{6}-2\theta)\sin^2{\frac{3}{12}\pi}}{\sin^2(\frac{3}{12}\pi-\theta)}~d\theta \\
  &<& 32\int_{-\frac{\pi}{12}}^{\frac{\pi}{12}}\cos^{2}(2\theta)
\sin^{2}(\frac{\pi}{6}+2\theta)\sin^2(\frac{\pi}{6}-2\theta)~d\theta\\
 &=& G.
\end{eqnarray*}

$(iii)$ \begin{eqnarray*}
  K_3 &=& 16\int_{-\frac{\pi}{12}}^{\frac{\pi}{12}}\frac{\cos^{2}(2\theta)
\sin^{2}(\frac{\pi}{6}+2\theta)\sin^2(\frac{\pi}{6}-2\theta)\sin^2({\frac{5}{12}\pi})}{\sin^2(\frac{5}{12}\pi-\theta)}~d\theta \\
  &<& \frac{2+\sqrt{3}}{6}\cdot G \\
  &<& G.
\end{eqnarray*}

The proof of Theorem \ref{prop 1 of TXY} is now complete.

\hfill $\Box$

\vspace{3mm}

\section{\textbf{Focal submanifolds with $g=6$.}}

\subsection{\textbf{On the focal submanifolds $M_1$ and $M_2$ with $(g, m_1, m_2)=(6, 1, 1)$}.}
\quad

This subsection will be devoted to a proof of Theorem \ref{thm1 TXY
focal g=6} (1).

Firstly, as mentioned before, the focal submanifolds are both
minimal in unit spheres. It follows that $\lambda_1(M_i)\leqslant
\dim M_i=5, i=1, 2.$ Next, we will only prove
$\lambda_1(M_1)\geqslant 3$, as the proof for $M_2$ is verbatim with
obvious changes on index ranges.

Recall the Dorfmeister-Neher theorem (\cite{DN}) which states that
the isoparametric hypersurface in $S^7(1)$ with $(g, m_1, m_2)=(6,
1, 1)$ is homogeneous. Further, as asserted by \cite{MO}, a
homogeneous hypersurface in $S^7(1)$ with $g=6$ is the inverse image
of the Cartan hypersurface in $S^4(1)$ with $g=3$ under the Hopf
fiberation (for the eigenvalues of Cartan hypersurfaces, see
\cite{Sol2}); this correspondence exists between focal submanifolds
of each hypersurface. Thus under the adjustment of the radius, we
get the following Riemannian submersion with totally geodesic
fibers:
\begin{eqnarray}\label{submersion g=6 m=1}
S^3(1)\hookrightarrow M_1 &\subset& S^7(1) \nonumber\\
    \downarrow~~~~ &&~~~~~~~ \downarrow \\
  S^2(\frac{\sqrt{3}}{2})/\mathbb{Z}_2 &\subset& S^4(\frac{1}{2}) \nonumber
\end{eqnarray}
where $S^2(\frac{\sqrt{3}}{2})/\mathbb{Z}_2 \subset
S^4(\frac{1}{2})$ is the Veronese embedding of the real projective
plane of constant Gaussian curvature $\frac{4}{3}$ into Euclidean
sphere of radius $\frac{1}{2}$.

Next, let us recall some background for the Laplacian of a
Riemannian submersion $\pi$ with totally geodesic fibers: $F
\hookrightarrow M \xrightarrow{\pi} B$. We denote the Laplacian of
$M$ by $\Delta^M$. At any point $m\in M$, the vertical Laplacian
$\Delta_v$ is defined to be
$$(\Delta_vf)(m)=((\Delta^{F_m})(f|_{F_m}))(m),$$ where
$F_m=\pi^{-1}(\pi(m))$ is the fiber of $\pi$ through $m$ and
$\Delta^{F_m}$ the Laplace operator of the metric induced by $M$ on
$F_m$. The horizontal Laplacian is the difference operator
$$\Delta_h=\Delta^M-\Delta_v.$$

According to Theorem 3.6 in \cite{BB}, the Hilbert space $L^2(M)$
admits a Hilbert basis consisting of simultaneous eigenfunctions for
$\Delta^M$ and $\Delta_v$. Then we can find a function $\phi$
satisfying:
\begin{equation*}
\left\{ \begin{array}{ll}
\Delta^{M_1}\phi=\lambda_1(M_1)\phi\\
~\Delta_v\phi~~=~b\phi.
\end{array}\right.
\end{equation*}
Since $\Delta_h\phi=(\lambda_1(M_1)-b)\phi$ and $\Delta_h$ is a
non-negative operator, we have
\begin{equation}\label{delta v leq 5}
b\leqslant \lambda_1(M_1)\leqslant 5.
\end{equation}
On the other hand, concerning the relation $Spec(\Delta_v)\subset
Spec(S^3(1))=\{0, 3, 8, ...\}$, we claim that $b\geqslant 3$. Otherwise,
suppose $b=0$, then $\phi$ is the composition of the fiberation
projection with an eigenfunction on the base space, such that
$$\lambda_1(M_1)\geqslant\lambda_1(S^2(\frac{\sqrt{3}}{2})/\mathbb{Z}_2)=8>5,$$
contradicting (\ref{delta v leq 5}).

Therefore, we arrive at
\begin{equation}\label{3 leq lambda_1 leq 5}
3\leqslant \lambda_1(M_1) \leqslant 5.
\end{equation}

\vspace{3mm}

\subsection{\textbf{The first eigenvalue of the focal submanifold $M_2$ with $(g, m_1, m_2)=(6,
2, 2)$}.}\label{4.2} \quad

Firstly, for sufficiently small $\varepsilon> 0$, we set
$$M_2(\varepsilon):=S^{n+1}(1)- B_{\varepsilon} (M_1)=\bigcup_{\theta\in[0, \frac{\pi}{6}-\varepsilon]}M_{\theta}$$
where $B_{\varepsilon} (M_1)=\{x\in S^{n+1}(1)~|~dist (x, M_1)<\varepsilon\}$, $M_{\theta}$ is the
isoparametric hypersurface with an oriented distance $\theta$ from $M_2$. Notice that the notation
$M_{\theta}$ here is different from that we used before.

Given $\theta\in (0,\frac{\pi}{6}-\varepsilon]$, let
$\{e_{\alpha,i}~ \boldsymbol{|} ~ i=1, 2,~ \alpha=1,..,6,~
e_{\alpha, i}\in E_{\alpha}\}$ be a local orthonormal frame field on
$M_{\theta}$ and $\xi$ be the unit normal field of $M_{\theta}$
towards $M_2$. After a parallel translation along the normal
geodesic from any point $x\in M_{\theta}$ to the point
$p=\phi_{\theta}(x)\in M_2$, (where $\phi_{\theta}:
M_{\theta}\rightarrow M_2$ is the focal map), the image of $\xi$ is
normal to the focal submanifold $M_2$ at $p$, which will still be
denoted by $\xi$; $e_{1,i}$ $(i=1, 2)$ turn out to be normal vectors
on $M_2$, which we will denote by $\widetilde{e}_{1,i},$ while the
others are still tangent vectors on $M_2$, which we will denote by
$\{\widetilde{e}_{2,i}, \widetilde{e}_{3,i}, \widetilde{e}_{4,i},
\widetilde{e}_{5,i}, \widetilde{e}_{6,i}\}$. They are determined by
$x$.

For any $X\in T_xM_{\theta}$, we can decompose it as
$X=X_1+\cdots+X_6\in E_1\oplus \cdots\oplus E_6$. Identify the
principal distribution $E_{\alpha}(x)$ ($\alpha=2,\cdots,6$, $x\in
M_{\theta}$) with its parallel translation at $p=\phi_{\theta}(x)\in
M_2$. The shape operator $A_{\xi}$ at $p$ is given in terms of its
eigenvectors $\widetilde{X}_{\alpha}$ (the parallel translation of
$X_{\alpha}, \alpha=2,\cdots,6)$ by (\emph{cf.} \cite{Mun})
\begin{eqnarray}\label{tilde X}
&&A_{\xi}\widetilde{X}_2=\cot(\theta_2-\theta_1)\widetilde{X}_2=\sqrt{3}\widetilde{X}_2, \quad A_{\xi}\widetilde{X}_3=\cot(\theta_3-\theta_1)\widetilde{X}_3=\frac{\sqrt{3}}{3}\widetilde{X}_3,\nonumber\\
&&A_{\xi}\widetilde{X}_4=\cot(\theta_4-\theta_1)\widetilde{X}_4=0, \qquad\quad A_{\xi}\widetilde{X}_5=\cot(\theta_5-\theta_1)\widetilde{X}_5=-\frac{\sqrt{3}}{3}\widetilde{X}_5,\\
&&A_{\xi}\widetilde{X}_6=\cot(\theta_6-\theta_1)\widetilde{X}_6=-\sqrt{3}\widetilde{X}_6.\nonumber
\end{eqnarray}
Namely, $\widetilde{X}_2, \widetilde{X}_3, \widetilde{X}_4, \widetilde{X}_5, \widetilde{X}_6$ belong to
the eigenspaces $E(\sqrt{3}), E(\frac{\sqrt{3}}{3}), E(0), E(-\frac{\sqrt{3}}{3}), E(-\sqrt{3})$ of $A_{\xi}$, respectively.

On the other hand, for any point $p\in M_2$, at a point $x\in
{\phi_{\theta}}^{-1}(p)$, the first principal distribution $E_1(x)$
is projected to be $0$ under ${(\phi_{\theta})}_{\ast}$; for the
others, we have
\begin{eqnarray*}
{(\phi_{\theta})}_{\ast}e_{\alpha,
i}&=&\frac{\sin(\theta_{\alpha}-\theta)}{\sin\theta_{\alpha}}\widetilde{e}_{\alpha,i}
=\frac{\sin\frac{\alpha-1}{6}\pi}{\sin(\frac{\alpha-1}{6}\pi+\theta)}\widetilde{e}_{\alpha,i}\\
&:=&\widetilde{k}_{\alpha-1}\widetilde{e}_{\alpha,i},~~\qquad\qquad~i=1, 2,~\alpha=2,\cdots,6.
\end{eqnarray*}
Denote by $\{\theta_{\alpha,i}~|~\alpha=1,\cdots,6,~ i=1, 2\}$ the
dual frame of $e_{\alpha,i}$. We then conclude that (up to a sign)
\begin{equation}\label{volume Mtheta}
dM_{\theta}=\prod_{j=1}^{2}\prod_{\alpha=2}^6\theta_{\alpha,j}\wedge
\prod_{i=1}^{2}\theta_{1,i}
=\frac{1}{(\widetilde{k}_1\cdots\widetilde{k}_5)^{2}}\phi_{\theta}^{\ast}(dM_2)\wedge
\prod_{i=1}^{2}\theta_{1,i}.
\end{equation}

Let $h$ be the same function as in last section. For sufficiently
small $\eta>0$, define $\widetilde{\psi}_{\eta}$ to be a nonnegative
smooth function on $[0, \frac{\pi}{2}-\eta]$ by
\begin{equation*}
\widetilde{\psi}_{\eta}(x):=\left\{ \begin{array}{ll}
1,\qquad\qquad x\in[0,\frac{\pi}{4}]\\
h(\frac{\frac{\pi}{2}-x}{\eta}),\quad x\in[\frac{\pi}{4}, \frac{\pi}{2}-\eta]
\end{array}\right.
\end{equation*}
Let $f_k$ $(k=0, 1, ...)$ be the $k$-th eigenfunctions on $M_2$ which are orthogonal to
each other with respect to the square integral inner product on $M_2$ and $L_{k+1} = Span\{f_0, f_1,..., f_k\}$.
Then any $\varphi\in L_{k+1}$ on $M_2$ can give rise to a function $\widetilde{\Phi}_{\varepsilon}: M_2(\varepsilon)\rightarrow \mathbb{R}$ by:
$$\widetilde{\Phi}_{\varepsilon}(x) = \widetilde{\psi}_{3\varepsilon}(3\theta)(\varphi\circ\phi_{\theta})(x).$$
Evidently, $\widetilde{\Phi}_{\varepsilon}$ is a smooth function on
$M_2(\varepsilon)$ satisfying the Dirichlet boundary condition and
square integrable on $M_2(\varepsilon)$.

As in last section, the calculation of $\|\nabla
\widetilde{\Phi}_{\varepsilon}\|^2_2$ is closely related to $|\nabla
\varphi(\phi_{\theta})|^2$.
%so we will firstly estimate $|\nabla \varphi(\phi_{\theta})|^2$.
According to (\ref{tilde X}), in the tangent space of $M_2$ at $p$,
we can decompose $\nabla \varphi$ as $\nabla \varphi =
Z_1+Z_2+Z_3+Z_4+Z_5\in E(\sqrt{3})\oplus E(\frac{\sqrt{3}}{3})\oplus
E(0)\oplus E(-\frac{\sqrt{3}}{3})\oplus E(-\sqrt{3})$. Subsequently,
\begin{equation}\label{nabla varphi rho}
\left\{ \begin{aligned}
\quad|\nabla \varphi|_p^2~~~&=|Z_1|^2+\cdots+|Z_5|^2\\
|\nabla
\varphi(\phi_{\theta})|_x^2&=\widetilde{k}^{2}_1|Z_1|^2+\cdots+\widetilde{k}^{2}_5|Z_5|^2
\end{aligned}\right.
\end{equation}

In the following, we intend to investigate the variation of $|\nabla
\varphi(\phi_{\theta})|^2$ along with the point $x$ in the fiber
sphere at $p$. For this purpose, we recall that each integral
submanifold of the curvature distributions corresponding to
$\cot\theta_j=\cot(\theta+\frac{j-1}{6}\pi)$ is a totally geodesic
submanifold in $M_{\theta}$ with constant sectional curvature $1 +
\cot^2\theta_j$ (\emph{cf}. for example, \cite{CCJ}). In our case,
we denote by $S^{2}(\sin\theta) \subset M_{\theta}$ the fiber sphere
at $p$. Then a similar calculation as in Section 3 leads to

\begin{eqnarray}\label{nabla Phi prime focal}
\lim_{\varepsilon\rightarrow 0}\|\nabla\widetilde{\Phi}_{\varepsilon}\|_2^2
&=&\lim_{\varepsilon\rightarrow 0}\int_{M_2(\varepsilon)}(\widetilde{\psi}_{3\varepsilon}(3\theta))^2
|\nabla (\varphi\circ \phi_{\theta})|^2 dM_2(\varepsilon)\\
&=&\int_{0}^{\frac{\pi}{6}}\Big(\int_{M_{\theta}}\frac{|\nabla
(\varphi\circ
\phi_{\theta})|^2}{{(\widetilde{k}_1}\cdots\widetilde{k}_5)^{2}}\phi_{\theta}^{\ast}(dM_2)dS^{2}
(\sin\theta)\Big)d\theta\nonumber
\end{eqnarray}

Given a point $p\in M_2\subset S^{13}(1)$, with respect to a
suitable tangent orthonormal basis $e_{\alpha}, e_{\bar{\alpha}}$
$(\alpha=1,...,5)$ of $T_pM_2$, as asserted by Miyaoka in
\cite{Miy}, the shape operators $A_{\xi}$, $A_{\zeta}$ and
$A_{\bar{\zeta}}$ with respect to the mutually orthogonal unit
normals: $\xi$ and two other unit normals, say $\zeta$ and
$\bar{\zeta}$, of $M_2$ are expressed respectively by diagonal
matrix
\begin{equation}\label{shape B(eta)}
A_{\xi}=\left(\begin{array}{ccccc}
\sqrt{3}I& & &  &  \\
& \frac{1}{\sqrt{3}}I&&&\\
&&0 & &\\
&&&-\frac{1}{\sqrt{3}}I&\\
&&&&-\sqrt{3}I
\end{array}\right)
\end{equation}
and symmetric matrices:
\begin{equation}\label{shape B(zeta)}
A_{\bar{\zeta}}=\left(\begin{array}{ccccc}
0  & -I & 0 & 0 & 0 \\
-I &  0 & 0 & \frac{2}{\sqrt{3}}I & 0  \\
0  &  0 & 0 & 0 & 0 \\
0  & \frac{2}{\sqrt{3}}I & 0 & 0 & -I\\
0  &  0 & 0 & -I & 0
\end{array}\right)\qquad
A_{\zeta}=\left(\begin{array}{ccccc}
0  &  J & 0 & 0 & 0 \\
-J &  0 & 0 & -\frac{2}{\sqrt{3}}J & 0  \\
0  &  0 & 0 & 0 & 0 \\
0  & \frac{2}{\sqrt{3}}J & 0 & 0 & J\\
0  &  0 & 0 & -J & 0
\end{array}\right)
\end{equation}
where \begin{equation}
I=\left(\begin{array}{cc}
1 & 0\\
0 & 1
\end{array}\right) \qquad
J=\left(\begin{array}{cc}
0 & -1\\
1 & 0
\end{array}\right).
\end{equation}

As a crucial step in our calculation, we set $\xi(t,s)=: \cos t~\xi
+\sin t \cos s ~\bar{\zeta} + \sin t\sin s ~\zeta$ ($0< t<
\pi,~0\leqslant s\leqslant 2\pi$) and the corresponding shape
operator $A(t, s)=: A_{\xi(t, s)}$, thus
\begin{equation}
A(t, s)=\left(\begin{array}{ccccc}
\sqrt{3}\cos t~ I & -\sin t~e^{-is} & 0 & 0 & 0\\
-\sin t~e^{is} & \frac{1}{\sqrt{3}}\cos t~I & 0 & \frac{2}{\sqrt{3}}\sin t~e^{-is} & 0\\
0& 0 & 0 & 0 & 0\\
0 & \frac{2}{\sqrt{3}}\sin t ~e^{-is} & 0 & -\frac{1}{\sqrt{3}}\cos t~I & -\sin t~e^{-is}\\
0 & 0 & 0 & -\sin t~e^{is} & -\sqrt{3}\cos t~ I
\end{array}\right),
\end{equation}
where $e^{is}$ is a matrix defined by $e^{is}=:\cos s ~I+\sin s ~J$.
The eigenvalues of $A(t, s)$ are still $\sqrt{3},
\frac{\sqrt{3}}{3}, 0, -\frac{\sqrt{3}}{3}$ and $-\sqrt{3}$, while
the corresponding eigenspaces of $A(t, s)$ are spanned by
eigenvectors as follows:

$E(\sqrt{3})=\mathrm{Span}\{\varepsilon_1, \varepsilon_{\bar{1}}\}$
with
\begin{eqnarray*}
 \varepsilon_1&=& \frac{1}{2\sqrt{2(1-\cos t)}}\Big( \sin t(1+\cos t)(\cos 2s~e_1-\sin 2s~e_{\bar{1}}) -\sqrt{3}\sin^2t(\cos s~e_2-\sin s~e_{\bar{2}})\\
&&\qquad\qquad\qquad\quad-\sqrt{3}\sin t(1-\cos t)e_4+(1-\cos t)^2(\cos s~ e_5+\sin s~e_{\bar{5}})\Big), \\
 \varepsilon_{\bar{1}}&=& \frac{1}{2\sqrt{2(1-\cos t)}}\Big( \sin t(1+\cos t)(\sin 2s~e_1+\cos 2s~e_{\bar{1}}) -\sqrt{3}\sin^2t(\sin s~e_2+\cos s~e_{\bar{2}})\\
&&\qquad\qquad\qquad\quad-\sqrt{3}\sin t(1-\cos t)e_{\bar{4}}+(1-\cos t)^2(-\sin s~ e_5+\cos s~e_{\bar{5}})\Big),
\end{eqnarray*}

$E(\frac{1}{\sqrt{3}})=\mathrm{Span}\{\varepsilon_2,
\varepsilon_{\bar{2}}\}$ with
\begin{eqnarray*}
 \varepsilon_2&=& \frac{1}{2\sqrt{2(1+\cos t)}}\Big( -\sqrt{3}\sin t(1+\cos t)(\cos 2s~e_1-\sin 2s~e_{\bar{1}}) \\
&&\qquad\qquad\qquad\quad+(1+\cos t)(1-3\cos t)(\cos s~e_2-\sin s~e_{\bar{2}})\\
&&\qquad\qquad\qquad\quad+\sin t(1+3\cos t)e_4+\sqrt{3}\sin^2t(\cos s~ e_5+\sin s~e_{\bar{5}})\Big), \\
 \varepsilon_{\bar{2}}&=& \frac{1}{2\sqrt{2(1+\cos t)}}\Big( -\sqrt{3}\sin t(1+\cos t)(\sin 2s~e_1+\cos 2s~e_{\bar{1}}) \\
&&\qquad\qquad\qquad\quad+(1+\cos t)(1-3\cos t)(\sin s~e_2+\cos s~e_{\bar{2}})\\
&&\qquad\qquad\qquad\quad+\sin t(1+3\cos t)e_{\bar{4}}+\sqrt{3}\sin^2t(-\sin s~ e_5+\cos s~e_{\bar{5}})\Big),
\end{eqnarray*}

$E(0)=\mathrm{Spann}\{\varepsilon_3, \varepsilon_{\bar{3}}\}$ with
$$\varepsilon_3=e_3, \quad \varepsilon_{\bar{3}}=e_{\bar{3}},$$

$E(-\frac{1}{\sqrt{3}})=\mathrm{Span}\{\varepsilon_4,
\varepsilon_{\bar{4}}\}$ with
\begin{eqnarray*}
 \varepsilon_4&=& \frac{1}{2\sqrt{2(1-\cos t)}}\Big( \sqrt{3}\sin t(1-\cos t)(\cos 2s~e_1-\sin 2s~e_{\bar{1}}) \\
&&\qquad\qquad\qquad\quad+(1-\cos t)(1+3\cos t)(\cos s~e_2-\sin s~e_{\bar{2}})\\
&&\qquad\qquad\qquad\quad+\sin t(1-3\cos t)e_4+\sqrt{3}\sin^2t(\cos s~ e_5+\sin s~e_{\bar{5}})\Big), \\
 \varepsilon_{\bar{4}}&=& \frac{1}{2\sqrt{2(1-\cos t)}}\Big( \sqrt{3}\sin t(1-\cos t)(\sin 2s~e_1+\cos 2s~e_{\bar{1}}) \\
&&\qquad\qquad\qquad\quad+(1-\cos t)(1+3\cos t)(\sin s~e_2+\cos s~e_{\bar{2}})\\
&&\qquad\qquad\qquad\quad+\sin t(1-3\cos t)e_{\bar{4}}+\sqrt{3}\sin^2t(-\sin s~ e_5+\cos s~e_{\bar{5}})\Big),
\end{eqnarray*}

$E(-\sqrt{3})=\mathrm{Span}\{\varepsilon_5, \varepsilon_{\bar{5}}\}$
with
\begin{eqnarray*}
 \varepsilon_5&=& \frac{1}{2\sqrt{2(1+\cos t)}}\Big( -\sin t(1-\cos t)(\cos 2s~e_1-\sin 2s~e_{\bar{1}}) -\sqrt{3}\sin^2t(\cos s~e_2-\sin s~e_{\bar{2}})\\
&&\qquad\qquad\qquad\quad+\sqrt{3}\sin t(1+\cos t)e_4+(1+\cos t)^2(\cos s~ e_5+\sin s~e_{\bar{5}})\Big), \\
 \varepsilon_{\bar{5}}&=& \frac{1}{2\sqrt{2(1+\cos t)}}\Big( -\sin t(1-\cos t)(\sin 2s~e_1+\cos 2s~e_{\bar{1}}) -\sqrt{3}\sin^2t(\sin s~e_2+\cos s~e_{\bar{2}})\\
&&\qquad\qquad\qquad\quad+\sqrt{3}\sin t(1+\cos
t)e_{\bar{4}}+(1+\cos t)^2(-\sin s~ e_5+\cos s~e_{\bar{5}})\Big).
\end{eqnarray*}

Now express $\nabla \varphi$ as
$$\nabla \varphi=\sum_{\alpha}(a_{\alpha}e_{\alpha}+a_{\bar{\alpha}}e_{\bar{\alpha}})
=\sum_{\alpha}(b_{\alpha}\varepsilon_{\alpha}+b_{\bar{\alpha}}\varepsilon_{\bar{\alpha}}),$$
where $a_{\alpha}=e_{\alpha}(\varphi)$,
$a_{\bar{\alpha}}=e_{\bar{\alpha}}(\varphi)$,
$b_{\alpha}=\varepsilon_{\alpha}(\varphi)$,
$b_{\bar{\alpha}}=\varepsilon_{\bar{\alpha}}(\varphi)$. It follows
that
\begin{eqnarray*}
%|\nabla \varphi|^2\quad &=& \sum_{\alpha=1}^5(a^2_{\alpha}+a^2_{\bar{\alpha}}) \\
|\nabla \varphi(\phi_{\theta})|^2 &=&
\sum_{\alpha=1}^5\widetilde{k}_{\alpha}^2(b^2_{\alpha}+b^2_{\bar{\alpha}}).
\end{eqnarray*}

Further, a direct calculation leads to
\begin{eqnarray*}
b^2_{1}+b^2_{\bar{1}} &=& I(t, s)+\frac{1}{8(1-\cos t)}\Big( \sin^2t(1+\cos t)^2(a_1^2+a_{\bar{1}}^2)+3\sin^4t(a_2^2+a_{\bar{2}}^2)\\
&&\qquad\qquad\qquad\qquad\qquad+3\sin^2t(1-\cos t)^2(a_4^2+a_{\bar{4}}^2)+(1-\cos t)^4(a_5^2+a_{\bar{5}}^2)\Big),
\\
b^2_{2}+b^2_{\bar{2}} &=& II(t, s)+\frac{1}{8(1+\cos t)}\Big( 3\sin^2t(1+\cos t)^2(a_1^2+a_{\bar{1}}^2)+(1+\cos t)^2(1-3\cos t)^2(a_2^2+a_{\bar{2}}^2)\\
&&\qquad\qquad\qquad\qquad\qquad+\sin^2t(1+3\cos t)^2(a_4^2+a_{\bar{4}}^2)+3\sin^4t(a_5^2+a_{\bar{5}}^2)\Big),
\\
b^2_{4}+b^2_{\bar{4}} &=&IV(t, s)+\frac{1}{8(1-\cos t)}\Big( 3\sin^2t(1-\cos t)^2(a_1^2+a_{\bar{1}}^2)+(1-\cos t)^2(1+3\cos t)^2(a_2^2+a_{\bar{2}}^2)\\
&&\qquad\qquad\qquad\qquad\qquad+\sin^2t(1-3\cos t)^2(a_4^2+a_{\bar{4}}^2)+3\sin^4t(a_5^2+a_{\bar{5}}^2)\Big),
\\
b^2_{5}+b^2_{\bar{5}} &=& V(t, s)+\frac{1}{8(1+\cos t)}\Big( \sin^2t(1-\cos t)^2(a_1^2+a_{\bar{1}}^2)+3\sin^4t(a_2^2+a_{\bar{2}}^2)\\
&&\qquad\qquad\qquad\qquad\qquad+3\sin^2t(1+\cos t)^2(a_4^2+a_{\bar{4}}^2)+(1+\cos t)^4(a_5^2+a_{\bar{5}}^2)\Big),
\end{eqnarray*}
where $I(t, s)$, $II(t, s)$, $IV(t, s)$ and $V(t, s)$ are those
items containing linear combinations of $\cos 2s\cos s$, $\cos
2s\sin s$, $\sin 2s\cos s$ and $\sin 2s\sin s$, whose integrals over
$s\in[0, 2\pi]$ vanish.

Transform $\int_{0}^{\pi}|\nabla\varphi\circ\phi_{\theta}|^2\sin
t~dt$ into
$$\int_{0}^{\pi}|\nabla\varphi\circ\phi_{\theta}|^2\sin
t~dt=(a_1^2+a_{\bar{1}}^2)A_1+(a_2^2+a_{\bar{2}}^2)A_2+(a_3^2+a_{\bar{3}}^2)A_3+(a_4^2+a_{\bar{4}}^2)A_4+(a_5^2+a_{\bar{5}}^2)A_5,$$
for some $A_1,\cdots,A_5.$ It is not difficult to find that
$$A_1=A_2=A_4=A_5=\frac{1}{2}(\widetilde{k}_1^2+\widetilde{k}_2^2+\widetilde{k}_4^2+\widetilde{k}_5^2).$$
Then we finally arrive at an estimate of
$\|\nabla\widetilde{\Phi_{\varepsilon}}\|_2^2$ in (\ref{nabla Phi
prime focal}):
\begin{eqnarray*}
\lim_{\varepsilon\rightarrow0}\|\nabla\widetilde{\Phi_{\varepsilon}}\|_2^2
&=&
\int_{0}^{\frac{\pi}{6}}\int_{M_2}\frac{\sin^2\theta}{(\widetilde{k}_1\cdots\widetilde{k}_5)^2}
\cdot2\pi\cdot\int_0^{\pi}|\nabla\varphi\circ\phi_{\theta}|^2\sin
t~dt dM_2 d\theta
\\
&=&
\int_{0}^{\frac{\pi}{6}}\int_{M_2}\frac{\sin^2\theta}{(\widetilde{k}_1\cdots\widetilde{k}_5)^2}
\cdot2\pi\cdot\Big(\frac{1}{2}(\widetilde{k}_1^2+\widetilde{k}_2^2+\widetilde{k}_4^2+\widetilde{k}_5^2)
(a_1^2+a_{\bar{1}}^2+a_2^2+a_{\bar{2}}^2\\
&&\qquad\qquad\qquad\qquad\qquad\qquad+a_4^2+a_{\bar{4}}^2+a_5^2+a_{\bar{5}}^2)+2\widetilde{k}_3^2(a_3^2+a_{\bar{3}}^2) \Big)dM_2d\theta
\\
&<&
\int_{0}^{\frac{\pi}{6}}\int_{M_2}\frac{\sin^2\theta}{(\widetilde{k}_1\cdots\widetilde{k}_5)^2}
\cdot2\pi\cdot\Big(\frac{1}{2}(\widetilde{k}_1^2+\widetilde{k}_2^2+\widetilde{k}_4^2+\widetilde{k}_5^2)
\sum_{\alpha=1}^5(a_{\alpha}^2+a_{\bar{\alpha}}^2) \Big)dM_2d\theta
\\
&=& \int_{0}^{\frac{\pi}{6}} \frac{4}{9}\sin^2\theta\cos^2\theta(\frac{1}{4}-\cos^22\theta)^2\Big(\frac{2-\cos 2\theta}{(\frac{1}{2}-\cos 2\theta)^2}
+\frac{3(2+\cos2\theta)}{(\frac{1}{2}+\cos 2\theta)^2}\Big) d\theta\cdot2\pi\cdot\|\nabla\varphi\|_2^2
\\
&=& (\frac{\pi}{18}-\frac{13\sqrt{3}}{240})\pi\cdot\|\nabla\varphi\|_2^2
\end{eqnarray*}
Combining with
\begin{eqnarray*}
\lim_{\varepsilon\rightarrow0}\|\widetilde{\Phi_{\varepsilon}}\|_2^2
&=&
\int_{0}^{\frac{\pi}{6}}\int_{M_2}\frac{\sin^2\theta}{(\widetilde{k}_1\cdots\widetilde{k}_5)^2}
\int_{M_2}\int_{S^{2}(\sin\theta)} \varphi(\phi_{\theta})^2 ~dS^{2}dM_2d\theta\nonumber\\
&=& 4\pi\cdot\|\varphi\|_2^2\cdot\int_0^{\frac{\pi}{6}}\frac{\sin^2\theta}{(\widetilde{k}_1\cdots\widetilde{k}_5)^2} ~d\theta \\
&=& \frac{16\pi}{9}\cdot\|\varphi\|_2^2\cdot \int_0^{\frac{\pi}{6}}\sin^22\theta(\frac{1}{4}-\cos^22\theta)^2~d\theta\\
&=& \frac{\pi^2}{108}\cdot\|\varphi\|_2^2
\end{eqnarray*}
we conclude that
\begin{eqnarray*}
\lim_{\varepsilon\rightarrow0}\frac{\|\nabla\widetilde{\Phi_{\varepsilon}}\|_2^2
}{\|\widetilde{\Phi_{\varepsilon}}\|_2^2} &<&
(6-\frac{117\sqrt{3}}{20\pi})\cdot\frac{\|\nabla\varphi\|_2^2}{\|\varphi\|_2^2}.
%&<&\frac{14}{5}~\frac{\|\nabla\varphi\|_2^2}{\|\varphi\|_2^2}.
\end{eqnarray*}
Similarly as the arguments in last section, we derive that
\begin{equation}\label{frac Phi prime}
\lambda_k(S^{13}(1))\leqslant
(6-\frac{117\sqrt{3}}{20\pi})~\lambda_k(M_2)<\frac{14}{5}~\lambda_k(M_2),
%\frac{(\frac{\pi}{18}-\frac{13\sqrt{3}}{240})}{\frac{\pi}{108}}~\lambda_k(M_2)<\frac{14}{5}~\lambda_k(M_2),
%\approx 2.774726<\frac{14}{5},
\end{equation}
as $\displaystyle
(6-\frac{117\sqrt{3}}{20\pi})\approx
2.774726.$ Taking $k=15$, the inequality turns to
$$\lambda_{15}(M_2)>10.$$

At last, recalling Lemma 3.1 in \cite{TY} which yields that the
dimension $10$ of $M_2$ is an eigenvalue of $M_2$ with multiplicity
at least $14$, we arrive at
$$\lambda_1(M_2)=\dim M_2=10~ \mathrm{with~ multiplicity}~ 14,$$
as required. The proof of Theorem \ref{thm1 TXY focal g=6} (ii) for
$M_2$ is now complete.

\vspace{3mm}

\subsection{\textbf{On the focal submanifold $M_1$ with $(g, m_1, m_2)=(6, 2, 2)$}.}
\quad

In this subsection, we still use the previous method to define
similar neighborhood $M_1(\varepsilon)$ of $M_1$ and the test
function $\widetilde{\Phi_{\varepsilon}}$. In the following, we will
just list the difference in the crucial step.

Given a point $ p\in M_1\subset S^{13}(1)$, with respect to a suitable
tangent orthonormal basis $e_{\alpha}, e_{\bar{\alpha}}$
$(\alpha=1,...,5)$ of $T_pM_1$, as asserted by Miyaoka in \cite{Miy}, the shape
operators $A_{\xi}$, $A_{\zeta}$ and $A_{\bar{\zeta}}$ with respect
to the mutually orthogonal unit normals: $\xi$ and two other unit normals, say $\zeta$ and
$\bar{\zeta}$, of $M_2$ are expressed respectively by symmetric
matrices:
\begin{equation}\label{shape B(eta)}
A_{\xi}=\left(\begin{array}{ccccc}
\sqrt{3}I& & &  &  \\
& \frac{1}{\sqrt{3}}I&&&\\
&&0 & &\\
&&&-\frac{1}{\sqrt{3}}I&\\
&&&&-\sqrt{3}I
\end{array}\right)
\end{equation}

\begin{equation}\label{shape B(zeta)}
A_{\bar{\zeta}}=\left(\begin{array}{ccccc}
0  & 0 & 0 & 0 & \sqrt{3}I \\
0 &  0 & 0 & \frac{1}{\sqrt{3}}I & 0  \\
0  &  0 & 0 & 0 & 0 \\
0  & -\frac{1}{\sqrt{3}}I & 0 & 0 & 0\\
-\sqrt{3}I &  0 & 0 & 0 & 0
\end{array}\right)\qquad
A_{\zeta}=\left(\begin{array}{ccccc}
0  &  0 & 0 & 0 & \sqrt{3}J \\
0  &  0 & 0 & \frac{1}{\sqrt{3}}J & 0  \\
0  &  0 & 0 & 0 & 0 \\
0  & -\frac{1}{\sqrt{3}}J & 0 & 0 & 0\\
-\sqrt{3}J  &  0 & 0 & 0 & 0
\end{array}\right)
\end{equation}
where \begin{equation}
I=\left(\begin{array}{cc}
1 & 0\\
0 & 1
\end{array}\right), \qquad
J=\left(\begin{array}{cc}
0 & -1\\
1 & 0
\end{array}\right).
\end{equation}

For the unit normal vector $\xi(t,s)=: \cos t~\xi +\sin t \cos s
~\zeta + \sin t\sin s ~\bar{\zeta}$ ($0<t<\pi,~0\leqslant s\leqslant
2\pi$), the corresponding shape operator $A(t, s)=: A_{\xi(t, s)}$
is given by
\begin{equation}
A(t, s)=\left(\begin{array}{ccccc}
\sqrt{3}\cos t~ I & 0 & 0 & 0 & \sqrt{3}\sin t~e^{i(\frac{\pi}{2}-s)}\\
0 & \frac{1}{\sqrt{3}}\cos t~I & 0 & \frac{1}{\sqrt{3}}\sin t~e^{i(\frac{\pi}{2}-s)} & 0\\
0& 0 & 0 & 0 & 0\\
0 & \frac{1}{\sqrt{3}}\sin t~e^{-i(\frac{\pi}{2}-s)} & 0 & -\frac{1}{\sqrt{3}}\cos t~I & 0\\
\sqrt{3}\sin t~e^{-i(\frac{\pi}{2}-s)} & 0 & 0 & 0 & -\sqrt{3}\cos t~ I
\end{array}\right).
\end{equation}
The eigenvalues of $A(t, s)$ are still $\sqrt{3},
\frac{\sqrt{3}}{3}, 0, -\frac{\sqrt{3}}{3}$ and $-\sqrt{3}$, while
the corresponding eigenspaces of $A(t, s)$ are spanned by
eigenvectors as follows:

$E(\sqrt{3})=\mathrm{Span}\{\varepsilon_1, \varepsilon_{\bar{1}}\}$ with
\begin{eqnarray*}
 \varepsilon_1&=& \frac{\sin t\sin s}{\sqrt{2(1-\cos t)}}e_1+\frac{\sin t\cos s}{\sqrt{2(1-\cos t)}}e_{\bar{1}}+\sqrt{\frac{1-\cos t}{2}}e_5\\
 \varepsilon_{\bar{1}}&=& -\frac{\sin t\cos s}{\sqrt{2(1-\cos t)}}e_1+\frac{\sin t\sin s}{\sqrt{2(1-\cos t)}}e_{\bar{1}}+\sqrt{\frac{1-\cos t}{2}}e_{\bar{5}},\\
\end{eqnarray*}

$E(\frac{1}{\sqrt{3}})=\mathrm{Span}\{\varepsilon_2, \varepsilon_{\bar{2}}\}$ with
\begin{eqnarray*}
 \varepsilon_2&=& \frac{\sin t\sin s}{\sqrt{2(1-\cos t)}}e_2+\frac{\sin t\cos s}{\sqrt{2(1-\cos t)}}e_{\bar{2}}+\sqrt{\frac{1-\cos t}{2}}e_4\\
 \varepsilon_{\bar{2}}&=& -\frac{\sin t\cos s}{\sqrt{2(1-\cos t)}}e_2+\frac{\sin t\sin s}{\sqrt{2(1-\cos t)}}e_{\bar{2}}+\sqrt{\frac{1-\cos t}{2}}e_{\bar{4}},\\
\end{eqnarray*}

$E(0)=\mathrm{Spann}\{\varepsilon_3, \varepsilon_{\bar{3}}\}$ with
$$\varepsilon_3=e_3, \quad \varepsilon_{\bar{3}}=e_{\bar{3}},$$

$E(-\frac{1}{\sqrt{3}})=\mathrm{Span}\{\varepsilon_4, \varepsilon_{\bar{4}}\}$ with
\begin{eqnarray*}
 \varepsilon_4&=& \frac{\sin t\sin s}{\sqrt{2(1+\cos t)}}e_2+\frac{\sin t\cos s}{\sqrt{2(1+\cos t)}}e_{\bar{2}}-\sqrt{\frac{1+\cos t}{2}}e_4\\
 \varepsilon_{\bar{4}}&=& -\frac{\sin t\cos s}{\sqrt{2(1+\cos t)}}e_2+\frac{\sin t\sin s}{\sqrt{2(1+\cos t)}}e_{\bar{2}}-\sqrt{\frac{1+\cos t}{2}}e_{\bar{4}},\\
\end{eqnarray*}

$E(-\sqrt{3})=\mathrm{Span}\{\varepsilon_5, \varepsilon_{\bar{5}}\}$ with
\begin{eqnarray*}
 \varepsilon_5&=& \frac{\sin t\sin s}{\sqrt{2(1+\cos t)}}e_1+\frac{\sin t\cos s}{\sqrt{2(1+\cos t)}}e_{\bar{1}}-\sqrt{\frac{1+\cos t}{2}}e_5\\
 \varepsilon_{\bar{5}}&=& -\frac{\sin t\cos s}{\sqrt{2(1+\cos t)}}e_1+\frac{\sin t\sin s}{\sqrt{2(1+\cos t)}}e_{\bar{1}}-\sqrt{\frac{1+\cos t}{2}}e_{\bar{5}},\\
\end{eqnarray*}

In an analogous way with that in last subsection, we obtain
\begin{eqnarray*}
\lim_{\varepsilon\rightarrow0}\|\nabla\widetilde{\Phi_{\varepsilon}}\|_2^2
&=&
\int_{0}^{\frac{\pi}{6}}\int_{M_1}\frac{\sin^2\theta}{(\widetilde{k}_1\cdots\widetilde{k}_5)^2}
\cdot2\pi\cdot\int_0^{\pi}|\nabla\varphi\circ\phi_{\theta}|^2\sin
t~dt dM_1 d\theta
\\
&=& 4\pi\int_{M_1}\Big((\frac{\pi}{36}-\frac{11\sqrt{3}}{240})(a_3^2+a_{\bar{3}}^2)
+(\frac{\pi}{144}+\frac{11\sqrt{3}}{1920})(a_1^2+a_{\bar{1}}^2+a_5^2+a_{\bar{5}}^2)\\
&&\qquad\quad+(\frac{\pi}{48}-\frac{21\sqrt{3}}{640})(a_2^2+a_{\bar{2}}^2+a_4^2+a_{\bar{4}}^2)\Big)~dM_1\\
&<&4\pi\int_{0}^{\frac{\pi}{6}}\frac{\sin^2\theta}{(\widetilde{k}_1\cdots\widetilde{k}_5)^2}\cdot
\frac{\widetilde{k}_1^2+\widetilde{k}_5^2}{2}~d\theta\cdot\|\nabla\varphi\|_2^2\\
&=&(\frac{\pi^2}{36}+\frac{11\sqrt{3}}{480}\pi)\cdot\|\nabla\varphi\|_2^2
\end{eqnarray*}
Combining with $\displaystyle\lim_{\varepsilon\rightarrow0}\|\widetilde{\Phi_{\varepsilon}}\|_2^2= \frac{\pi^2}{108}\cdot\|\varphi\|_2^2$,
we eventually arrive at
$$\lambda_k(S^{13}(1))\leqslant\lambda_k(M_1)\cdot \frac{\frac{\pi^2}{36}+\frac{11\sqrt{3}}{480}\pi}{\frac{\pi^2}{108}}=\lambda_k(M_1)\cdot(3+\frac{99\sqrt{3}}{40\pi}),$$
as required. This completes the proof of Theorem \ref{thm1 TXY focal
g=6} (ii) for $M_1$.

\vspace{4mm}

\section{\textbf{Focal submanifolds with $g=4$.}}

We begin this section with a short review of the isoparametric
hypersurfaces of OT-FKM-type. For a symmetric Clifford system
$\{P_0,\cdots,P_m\}$ on $\mathbb{R}^{2l}$, \emph{i.e.}, $P_i$'s are
symmetric matrices satisfying $P_iP_j+P_jP_i=2\delta_{ij}I_{2l}$,
Ferus, Karcher and M\"{u}nzner (\cite{FKM}) constructed a polynomial
$F$ on $\mathbb{R}^{2l}$:
\begin{eqnarray}\label{FKM isop. poly.}
&&\qquad F:\quad \mathbb{R}^{2l}\rightarrow \mathbb{R}\nonumber\\
&&F(x) = |x|^4 - 2\displaystyle\sum_{i = 0}^{m}{\langle
P_ix,x\rangle^2}.
\end{eqnarray}

It turns out that each level hypersurface of $f=F|_{S^{2l-1}}$,
\emph{i.e.}, the preimage of some regular value of $f$, is an
isoparametric hypersurface with four distinct constant principal
curvatures. Choosing $\xi=\frac{\nabla f}{|\nabla f|}$, we have
$M_1=f^{-1}(1)$, $M_2=f^{-1}(-1)$, which have codimensions $m_1+1$
and $m_2+1$ in $S^{n+1}(1)$, respectively. The multiplicity pairs
$(m_1, m_2)=(m, l-m-1)$, provided $m>0$ and $l-m-1> 0$, where $l =
k\delta(m)$ $(k=1,2,3,...)$ and $\delta(m)$ is the dimension of an
irreducible module of the Clifford algebra $\mathcal{C}_{m-1}$ on
$\mathbb{R}^{l}$.

\subsection{\textbf{On the focal submanifold $M_1$ with $(g, m_1, m_2)=(4, 1, 1)$.}}\quad

As mentioned before, the isoparametric foliation with $(g, m_1,
m_2)=(4, 1, 1)$ can be expressed in the form of OT-FKM-type. In
fact, by an orthogonal transformation, we can always choose the
Clifford matrices $P_0$, $P_1$ to be
\begin{equation*}
P_0=\left(\begin{array}{cc}
  1 & 0 \\
  0 & -1
\end{array}\right)
,
P_1=\left(\begin{array}{cc}
  0 & 1 \\
  1 & 0
\end{array}\right).
\end{equation*}
Then the focal submanifold $M_1$ is expressed as
$$M_1=:\{(x, y) \in \mathbb{R}^3\times\mathbb{R}^3~|~\langle x, y\rangle=0, |x|=|y|=\frac{1}{\sqrt{2}}\}.$$

In order to investigate $M_1$, we define a two-fold covering as
follows, regarding $S^3$ as the group of unit vectors in
$\mathbb{H}$ of quaternions:
\begin{eqnarray*}
\sigma: S^3 &\rightarrow& M_1\subset \mathbb{R}^6 \\
a &\mapsto& \frac{1}{\sqrt{2}}(aj\bar{a}, ak\bar{a})
\end{eqnarray*}
where $i$, $j$, $k$ are basis elements satisfying $i^2=j^2=k^2=-1$
and $ij=k$.

Let us equip $S^3$ with the induced metric by $\sigma$. To be more
specific, at any point $a\in S^3$, we can choose a basis of $T_aS^3$
as $e_1=:ai$, $e_2=:aj$, $e_3=:ak$, whose images under the tangent
map $\sigma_{\ast}(X)=\frac{1}{\sqrt{2}}(aj\bar{X}+Xj\bar{a},
ak\bar{X}+Xk\bar{a})$ ($X\in T_aS^3$) are
\begin{eqnarray*}
&&\sigma_{\ast}(e_1)=(\sqrt{2}ak\bar{a}, -\sqrt{2}aj\bar{a}),\,\,
\sigma_{\ast}(e_2)=(0, \sqrt{2}ai\bar{a}),\,\,
\sigma_{\ast}(e_3)=(-\sqrt{2}ai\bar{a}, 0).
\end{eqnarray*}
Subsequently, the metric matrix is
%the metric matrix of $\Big(\langle e_p, e_q \rangle\Big)$ $(p, q=1, 2, 3)$ is easily
\begin{equation*}
\Big(\langle e_p, e_q \rangle\Big)=\left(\begin{array}{ccc}
4&0&0\\
0&2&0\\
0&0&2
\end{array}
\right).
\end{equation*}

Therefore, $S^3$ with the induced metric is a Berger sphere, say
$S^3_B$, and $M_1$ is isometric to the $\mathbb{Z}_2$-quotient
$S^3_B/\mathbb{Z}_2$ by identifying its antipodal points. Actually,
identifying $\mathbb{C}\times\mathbb{C}$ with $\mathbb{H}$ by $(z,
w)\rightarrow z+jw$, we have the Hopf fiberation $S^3_B\rightarrow
S^2$ defined by $(z, w)\mapsto (|z|^2-|w|^2, 2z\bar{w})$, which
gives rise to the following Riemannian submersion with totally
geodesic fibers and with the vertical space spanned by $e_1$:
\begin{eqnarray}\label{Berger sphere}
S^1(1)\hookrightarrow & S^3_B/\mathbb{Z}_2&\cong M_1\nonumber\\
&\downarrow&\\
&S^2(\frac{\sqrt{2}}{2})&\nonumber
\end{eqnarray}
Comparing with the Riemannian submersion with totally geodesic
fibers
\begin{eqnarray}\label{Rie sub S3}
S^1(\frac{\sqrt{2}}{2})\hookrightarrow & S^3(\sqrt{2})/\mathbb{Z}_2&\nonumber\\
&\downarrow\pi&\\
&S^2(\frac{\sqrt{2}}{2})&\nonumber
\end{eqnarray}
where $S^3(\sqrt{2})$ is the standard sphere with radius $\sqrt{2}$,
we can calculate the first eigenvalue of $M_1\cong
S^3_B/\mathbb{Z}_2$ in the following steps.

Firstly, given a Riemannian submersion $\pi: (M,g) \rightarrow B$
with totally geodesic fibers, for each $t>0$, there is a unique
Riemannian metric $g_t$ on $M$, such that for any $m\in M$,
\begin{enumerate}
  \item $g_t|_{V_mM\times H_mM}=0$;
  \item $g_t|_{V_mM}=t^2g|_{V_mM}$;
  \item $g_t|_{H_mM}=g|_{H_mM}$.
\end{enumerate}
We denote by $M_{g_t}$ the Riemannian manifold $(M, g_t)$ and by
$\Delta_t^{M}$ its Laplacian. It is clear that
$\Delta_t^{M}=t^{-2}\Delta_v+\Delta_h$ (\emph{cf}. \cite{BB}). Thus,
a common eigenfunction of $\Delta_v$ and $\Delta_h$ is an
eigenfunction of $\Delta_t^M$.

In contrast with our case, we see that
$M=S^3(\sqrt{2})/\mathbb{Z}_2$ and $M_{g_t}=S^3_B/\mathbb{Z}_2$ with
$t=\sqrt{2}$.

Secondly, denote the spectrum of the Riemannian manifold $M$ by
$\{(\mu_k, n_k)~|~ 0=\mu_0<\mu_1<\cdots<\mu_k<\cdots \uparrow
\infty; \mbox{$\mu_k$ is an eigenvalue, $n_k$ is the multiplicity of
$\mu_k$}\}.$ For the convenience, we list the well known spectrums
of $S^1(1), S^1(\frac{\sqrt{2}}{2}), S^2(\frac{\sqrt{2}}{2}),
S^3(\sqrt{2})/\mathbb{Z}_2$ in Table \ref{table 1}.

\begin{table}[htbp]
\caption{}\label{table 1}
\begin{tabular}{|c|c|c|c|c|}
  \hline
   $M$      & $(\mu_1,n_1)$ & $(\mu_2,n_2)$ & $(\mu_3,n_3)$ & $(\mu_k,n_k)(k>1)$ \\
  \hline
  $S^1(1)$ & $(1,2)$       & $(4,2)$ & $(9,2)$ & $(k^2,2)$  \\\hline
  $S^1(\frac{\sqrt{2}}{2})$ & $(2,2)$ & $(8,2)$ & $(18,2)$ & $(2k^2,2)$  \\\hline
  $S^2(\frac{\sqrt{2}}{2})$ & $(4,3)$ & $(12,5)$ & $(24,7)$ & $(2k(k+1),2k+1)$  \\\hline
  $S^3(\sqrt{2})/\mathbb{Z}_2$ & $(4,9)$ & $(12,25)$ & $(24,49)$ & $(2k(k+1),(2k+1)^2)$ \\
  \hline
\end{tabular}
\end{table}

Finally, let $\Delta_h$, $\Delta_v$ be the corresponding horizontal
and vertical Laplacians in the Riemannian submersion (\ref{Rie sub
S3}). From Theorem 3.6 in \cite{BB}, it follows that for any
$\lambda\in Spec(S^3(\sqrt{2})/\mathbb{Z}_2)$, there exist
nonnegative real numbers $b\in Spec(\Delta_h)$ and $\phi\in
Spec(\Delta_v)\subset Spec(S^1(\frac{\sqrt{2}}{2}))$, such that
$\lambda=b+\phi$. As discussed at the first step, we see that
$\bar{\lambda}:=b+\frac{1}{2}\phi\in Spec(S^3_B/\mathbb{Z}_2)$.
According to Table \ref{table 1}, there are only three cases to be
considered:

$(i)$ $\lambda=0$. Obviously, in this case $\bar{\lambda}=0$.

$(ii)$ $\lambda\geqslant 12$. We claim that $\bar{\lambda}>3$.
Suppose $\bar{\lambda}\leqslant 3$. Then the inequality
$\frac{1}{2}\phi\leqslant\bar{\lambda}\leqslant3$ implies that
$\phi=0$ or $2$. Hence $\bar{\lambda}\geqslant
b=\lambda-\phi\geqslant10$, which contradicts the assumption.

$(iii)$ $\lambda=4$. Clearly, $b, \phi \geqslant 0$ and $4=b+\phi$.
From Table $\ref{table 1}$, it follows that the possible values of
$\phi$ are only $0$ or $2$. Let $E_1$ be the eigenspace
corresponding to $\lambda=4$. Again by Theorem 3.6 in \cite{BB},
there exist linearly independent functions $f_1,\cdots,f_9$ such
that $E_1=\mathrm{Span}\{f_1,\cdots,f_9\}$ and
\[  \Delta_hf_k=b_kf_k, ~~\Delta_vf_k=\phi_kf_k, ~~b_k+\phi_k=4,\hspace{10pt}\mbox{for}~k=1,2,\cdots,9.       \]
Let $i$ be the non-negative integer such that $\phi_k=0$, for
$k\leqslant i$; $\phi_k=2$, for $k>i$. If $k\leqslant i$, the
corresponding function $f_k$ is induced from the base space. That
is, there exists some function $h_k$ such that $f_k=h_k\circ\pi,~
\Delta_Bh_k=4h_k,$ where $\Delta_B$ is the Laplacian on the base
manifold $S^2(\frac{\sqrt{2}}{2})$. Since the multiplicity of $4\in
Spec(S^2(\frac{\sqrt{2}}{2}))$ is $3$, it yields that $i=3$. Namely,
$\phi_k=2$ and $b_k=2$ for $k>3$. Subsequently,
$$ b_k+\frac{1}{2}\phi_k=3\in Spec(S^3_B/\mathbb{Z}_2). $$
Moreover, the space consisting of such functions has dimension $6$.

Putting all these facts together, we complete the proof of the first part in Theorem \ref{thm2 TXY focal g=4}.

\vspace{3mm}

\subsection{\textbf{On the homogeneous focal submanifold $M_1$ with $(g, m_1, m_2)=(4, 4, 3)$.}}
\quad

The last subsection will be devoted to calculating the first
eigenvalue of the focal submanifold $M_1$ with dimension $10$ and
$(g, m_1, m_2)=(4, 4, 3)$ of OT-FKM type in $S^{15}(1)$. We use
analogous method as that in Subsection \ref{4.2} to define
$M_1(\varepsilon)$ and $\widetilde{\Phi_{\varepsilon}}$. In the
following, we calculate
$\|\nabla\widetilde{\Phi}_{\varepsilon}\|_2^2$.

Firstly, let us make some notations. For any $x\in M_1$, denote
$\nabla\varphi|_x=:X\in T_xM_1$. To simplify the illustration, we
assume temporarily $|X|=1$. For any point $a=(a_0,\cdots,a_4)$ in
the unit sphere $S^4(1)$, let
$P_a=:\sum_{\beta=0}^4a_{\beta}P_{\beta}$ be an element in the
Clifford sphere $\Sigma=:\Sigma(P_0,\cdots, P_4)$ spanned by
$P_0,\cdots,P_4$. Denote $\xi_a=:P_ax$. Then its shape operator is
$A_{\xi_a}=\sum_{\beta=0}^4a_{\beta}A_{\xi_{\beta}}$.

Next, in virtue of \cite{FKM}, for any $a\in S^4(1)$, we can
decompose $X$ with respect to eigenspaces of $A_{\xi_a}$ into
$$X=Y_1+Y+Y_{-1}\in E_{1}(A_{\xi_a})\oplus E_0(A_{\xi_a})\oplus E_{-1}(A_{\xi_a}).$$
Recall that
$T^{\perp}_xM_1=\mathrm{Span}\{P_{\beta}x~|~\beta=0,\cdots,4\}$ and
$E_0(A_{\xi_a})=\mathbb{R}\{QP_ax~|~Q\in\Sigma, \langle Q,
P_a\rangle=0\}$. Thus if we choose $Q_j$ $(j=1,\cdots,4)$ in such a
way that they constitute with $P_a$ an orthonormal basis of
$\Sigma$, then $Y=\sum_{j=1}^4\langle X, Q_jP_ax\rangle Q_jP_ax$,
and hence
\begin{equation}\label{|Y|^2}
|Y|^2=\sum_{j=1}^4\langle X, Q_jP_ax\rangle ^2=\sum_{j=1}^4\langle
P_aX, Q_jx\rangle^2=|(P_aX)^{\perp}|^2.
\end{equation}
Therefore, combining with the formula $A_{\xi_a}X=-(P_aX)^T$,
%combining with the decomposition $P_aX=(P_aX)^T+(P_aX)^{\perp}$,
we get $|Y|^2=1-|A_{\xi_a}X|^2$.

On the other hand, notice that
$$|A_{\xi_a}X|^2=\sum_{\alpha, \beta=0}^4a_{\alpha}a_{\beta}\langle A_{\xi_{\alpha}}X, A_{\xi_{\beta}}X \rangle
=\sum_{\beta=0}^4a_{\beta}^2|A_{\xi_{\beta}}X|^2+T,$$ where $T$ is
the item consisting of the products $a_{\alpha}a_{\beta}\langle
A_{\xi_{\alpha}}X, A_{\xi_{\beta}}X \rangle~(\alpha\neq \beta)$,
whose integral on $S^4(1)$ vanishes since
$\int_{S^4(1)}a_{\alpha}a_{\beta}dv=0$ for $\alpha\neq \beta$. By
the decomposition
$$P_{\beta}X=(P_{\beta}X)^T+(P_{\beta}X)^{\perp}=(P_{\beta}X)^T+\sum_{\gamma=0}^4\langle P_{\beta}X, P_{\gamma}x\rangle P_{\gamma}x,$$
we obtain that $1=|A_{\xi_{\beta}}X|^2+\sum_{\gamma=0}^4\langle
P_{\beta}X, P_{\gamma}x\rangle^2$.
%From the relation $A_{\xi_a}X=\sum_{\beta=0}^4a_{\beta}A_{\xi_{\beta}}X$,
Therefore, the arguments above imply that
$$|Y|^2=1-\sum_{\beta=0}^4a_{\beta}^2|A_{\xi_{\beta}}X|^2-T=1-\sum_{\beta=0}^4a_{\beta}^2(1-\displaystyle\sum_{\gamma=0}^4\langle P_{\beta}X, P_{\gamma}x\rangle^2)-T=\sum_{\beta, \gamma=0}^4a_{\beta}^2\langle P_{\beta}X, P_{\gamma}x\rangle^2-T,$$
which leads to
\begin{eqnarray*}
\int_{S^4(1)}|Y|^2dv&=&\int_{S^4(1)}\sum_{\beta,
\gamma=0}^4a_{\beta}^2\langle P_{\beta}X, P_{\gamma}x\rangle^2dv
=\frac{1}{5}\mathrm{Vol(S^4(1))}\cdot\sum_{\beta, \gamma=0}^4\langle
X,
P_{\beta}P_{\gamma}x\rangle^2\\
&=&\frac{2}{5}\mathrm{Vol(S^4(1))}, \end{eqnarray*} where the last
equality is followed from a crucial assertion that
$$T_xM_1=\mathrm{Span}\{P_{\beta}P_{\gamma}x~|~ \beta,
\gamma=0,\cdots,4, ~~\beta<\gamma\}$$ which holds only for
homogeneous case with $(g, m_1, m_2)=(4, 4, 3)$(\emph{cf}.
Subsection 3.2.1 1), the case $q=2$ in \cite{QTY}). Subsequently, it
is easily seen that
$$\int_{S^4(1)}|Y_1|^2dv=\int_{S^4(1)}|Y_{-1}|^2dv=\frac{1}{2}\int_{S^4(1)}(|X|^2-|Y|^2)dv=\frac{3}{10}\mathrm{Vol(S^4(1))}.$$

In this way, we obtain that
\begin{eqnarray*}
\lim_{\varepsilon\rightarrow 0}\|\nabla\widetilde{\Phi}_{\varepsilon}\|_2^2 &=& \int_0^{\frac{\pi}{4}}\int_{M_{\theta}}|\nabla\varphi(\phi_{\theta})|^2~dM_{\theta}d\theta \\
 &=& \int_{0}^{\frac{\pi}{4}}\int_{M_1}\int_{S^4(\sin\theta)}\frac{1}{\widetilde{k}_1^3\widetilde{k}_2^4\widetilde{k}_3^3}
 (\widetilde{k}_1^2|Y_1|^2+\widetilde{k}_2^2|Y|^2+\widetilde{k}_3^2|Y_{-1}|^2)~dS^4(\sin\theta)dM_{1}d\theta \\
 &=& \int_{0}^{\frac{\pi}{4}}\int_{M_1}\mathrm{Vol(S^4(\sin\theta))}|\nabla\varphi|^2\Big(\frac{3}{10}(\frac{1}{(\cos \theta+\sin\theta)^2}+\frac{1}{(\cos\theta-\sin\theta)^2})\\
 &&+\frac{2}{5}\frac{1}{\cos^2\theta}\Big)\cdot\cos^32\theta\cos^4\theta ~dM_{1}d\theta \\
 &=& \|\nabla\varphi\|_2^2\cdot\mathrm{Vol(S^4(1))}\Big( \frac{17}{2400}-\frac{\pi}{1280}\Big).
\end{eqnarray*}
Combining with
\begin{eqnarray}\label{Phi prime}
\lim_{\varepsilon\rightarrow
0}\|\widetilde{\Phi}_{\varepsilon}\|_2^2
&=&\int_0^{\frac{\pi}{4}}\frac{1}{\widetilde{k}_1^{3}\widetilde{k}_2^{4}\widetilde{k}_3^{3}} \int_{M_1}\int_{S^{4}(\sin\theta)} \varphi(\phi_{\theta})^2 dS^{4}(\sin\theta)dM_1d\theta\nonumber\\
&=& \|\varphi\|_2^2\cdot\int_0^{\frac{\pi}{4}}\frac{1}{\widetilde{k}_1^{3}\widetilde{k}_2^{4}\widetilde{k}_3^{3}} Vol(S^{4}(\sin\theta))~d\theta \nonumber\\
&=&\|\varphi\|_2^2\cdot \mathrm{Vol(S^4(1))}\frac{1}{64}\cdot B(\frac{5}{2}, 2)\nonumber\\
&=&\|\varphi\|_2^2\cdot \mathrm{Vol(S^4(1))}\frac{1}{560},\nonumber
\end{eqnarray}
we conclude that
\begin{equation*}
\lim_{\varepsilon\rightarrow
0}\frac{\|\nabla\widetilde{\Phi}_{\varepsilon}\|_2^2}{\|\widetilde{\Phi}_{\varepsilon}\|_2^2}
%=\frac{ \|\nabla\varphi\|_2^2}{ \|\varphi\|_2^2}\cdot560\cdot(\frac{17}{2400}-\frac{\pi}{1280})
=\frac{ \|\nabla\varphi\|_2^2}{
\|\varphi\|_2^2}\cdot(\frac{119}{30}-\frac{7\pi}{16})
\end{equation*}

Analogously as in Section 4,
\begin{equation*}
\lambda_k(S^{15}(1))\leqslant
\lambda_k(M_1)\cdot(\frac{119}{30}-\frac{7\pi}{16}).
\end{equation*}
In particular, $$\lambda_{17}(M_1)\geqslant
\frac{32}{\frac{119}{30}-\frac{7\pi}{16}}>12.$$

Finally, combing with Lemma 3.1 in \cite{TY}, we conclude that
$$\lambda_1(M_1)=\dim M_1=10,~ \mathrm{with~ multiplicity}~ 16,$$
as required.

%%%%%%%%%%%%%%%%%%%%%%%%%%%%%%%%%%%%%%%%%%%%%%%%%%%%%%%%%%%%%%%%%%%%%%%%%%%%%%%%%%%%%%%%%%%%%%%%%%%%%%%%%%%%%%%%%%%%%%%%%%%%%%%%%%%%%%%%%

\begin{ack}
We would like to thank Professor Reiko Miyaoka sincerely for her
interest. We also would like to thank the anonymous referees for
valuable comments.
\end{ack}


\begin{thebibliography}{123}

\bibitem[Abr]{Abr}
U. Abresch, \emph{Isoparametric hypersurfaces with four or six distinct principal curvatures}, Math.
Ann. 264(1983), 283-302.

\bibitem[BB]{BB}
L. B\'{e}rard Bergery and J.-P. Bourguignon, \emph{Riemannian submersions with totally geodesic
fibers}, Illinois J. Math. \textbf{26} (1982), 181-¨C200.

\bibitem[Car1]{Car1}
E. Cartan, \emph{Sur des familles remarquables d'hypersurfaces isoparam\'{e}triques dans les
espaces sph\'{e}riques}, Math. Z. \textbf{45} (1939), 335--367.

\bibitem[Car2]{Car2}
E. Cartan, \emph{Sur quelque familles remarquables d'hypersurfaces}, C. R. Congr\`{e}s Math.
Li\`{e}ge, 1939, 30--41.

\bibitem[CCJ]{CCJ}
T. E. Cecil, Q.-S. Chi and G. R. Jensen, \emph{Isoparametric hypersurfaces with four
principal curvatures}, Ann. Math. \textbf{166} (2007), 1--76.

\bibitem[CF]{CF}
I.~ Chavel and E.~A.~ Feldman: \emph{Spectra of domains in compact manifolds},
J. Funct. Anal. \textbf{30} (1978), 198--222.


\bibitem[Chi]{Chi}
Q. S. Chi, \emph{Isoparametric hypersurfaces with four principal
curvatures, III}, J. Differential Geom. \textbf{94} (2013),
469--504.

\bibitem[CR]{CR}
T.~E. Cecil and P.~T. Ryan, \emph{Tight and taut immersions of
manifolds}, Research Notes in Math. \textbf{107}, Pitman, London, 1985.

\bibitem[DN]{DN}
J. Dorfmeister and E. Neher, \emph{Isoparametric hypersurfaces, case
$g = 6$, $m = 1$}, Comm. in Algebra \textbf{13} (1985), 2299--2368.

\bibitem[FKM]{FKM} D. Ferus, H. Karcher, and H.~F. M{\"u}nzner,
\emph{Cliffordalgebren und neue isoparametrische Hyperfl\"{a}chen},
Math. Z. \textbf{177} (1981), 479--502. For an English version, see arXiv: 1112.2780.

\bibitem[Imm]{Imm}
S. Immervoll, \emph{On the classification of isoparametric
hypersurfaces with four distinct principal curvatures in spheres},
Ann. Math. \textbf{168} (2008), 1011--1024.

\bibitem[Kot]{Kot}
M.~Kotani, \emph{The first eigenvalue of homogeneous minimal hypersurfaces in a unit sphere $S^{n+1}(1)$},
T\^{o}hoku Math. J. \textbf{37} (1985), 523--532.

\bibitem[Miy1]{Miy2}
R. Miyaoka, \emph{The Dorfmeister-Neher theorem on isoparametric
hypersurfaces}, Osaka J. Math. \textbf{46} (2009), 695--715.

\bibitem[Miy2]{Miy}
R. Miyaoka, \emph{Isoparametric hypersurfaces with $(g,m) = (6,
2)$}, Ann. Math. \textbf{177} (2013), 53-110.

\bibitem[MO]{MO}
R. Miyaoka and T. Ozawa, \emph{Construction of taut embeddings and
Cecil-Ryan conjecture}, in ``Geometry of manifolds (edited by
K.Shiohama)", Acad. Press, London, 1989, 181--189.

\bibitem[MOU]{MOU}
H.~Muto, Y.~Ohnita, and H.~Urakawa, \emph{Homogeneous minimal
hypersurfaces in the unit sphere and the first eigenvalue of the
Laplacian}, T\^{o}hoku Math. J. \textbf{36} (1984), 253--267.

\bibitem[M{\"u}n]{Mun}
H. F. M{\"u}nzner, \emph{Isoparametric hyperfl\"achen in sph\"aren, I
and II},~Math. Ann. \textbf{251} (1980), 57--71 and \textbf{256} (1981), 215--232.

\bibitem[Mut]{Mut}
H.~Muto, \emph{The first eigenvalue of the Laplacian of an isoparametric minimal hypersurface in a unit sphere},
Math. Z.\textbf{ 197} (1988), 531--549

\bibitem[Oza]{Oza}
S.~ Ozawa, \emph{Singular variation of domains and eigenvalues of the Laplacian},
Duke Math. J. \textbf{48} (1981), 767--778.

\bibitem[QTY]{QTY}
C. Qian, Z. Z. Tang and W. J. Yan, \emph{New examples of Willmore
submanifolds in the unit sphere via isoparametric functions, II},
Ann. Glob. Anal. Geom. \textbf{43}, (2013), 47-62.


\bibitem[Sol1]{Sol}
B.~Solomon, \emph{Quartic isoparametric hypersurfaces and quadratic forms}, Math. Ann. \textbf{293} (1992), 387--398.


\bibitem[Sol2]{Sol2}
B.~Solomon, \emph{The harmonic analysis of cubic isoparametric
minimal hypersurfaces I: Dimensions $3$ and $6$; II: Dimensions $12$
and $24$}, Amer. J. of Math. \textbf{112 }(1990), 151--203;
205--241.

\bibitem[TY]{TY}
Z.~Z.~Tang and W.~J.~Yan, \emph{Isoparametric foliation and Yau
conjecture on the first eigenvalue}, J. Differential Geom.
\textbf{94} (2013), 521--540.


\bibitem[Yau]{Yau}
S. T. Yau, Problem section, \emph{Seminar on differential geometry}, Ann. Math. Studies \textbf{102},
Princeton Univ. Press, 1982.


\end{thebibliography}
\end{document}